\documentclass{gtart}

\usepackage{epsf,amssymb}

\def\ifplaintex{\expandafter\ifx\csname documentclass\endcsname\relax}

\ifplaintex 
\else
\fi

\expandafter\ifx\csname beginpicture\endcsname\relax
\expandafter\ifx\csname documentclass\endcsname\relax
\input pictex \else
\input prepictex \input pictex \input postpictex \fi\fi

\def\gt{{\mathsurround=0pt\it $\cal G\mskip-2mu$eometry \&\ 
$\cal T\!\!$opology}}        

\def\gtp{{\mathsurround=0pt\it $\cal G\mskip-2mu$eometry \&\ 
$\cal T\!\!$opology $\cal P\!$ublications}}  

\def\volumenumber#1{\def\thevolumenumber{#1}}
\def\papernumber#1{\def\thepapernumber{#1}}
\def\volumeyear#1{\def\thevolumeyear{#1}}

\def\pagenumbers#1#2{\def\startpage{#1}\def\finishpage{#2}}
\def\published#1{\def\publishdate{#1}}
\def\proposed#1{\def\theproposer{#1}}
\def\seconded#1{\def\theseconders{#1}}
\def\received#1{\def\receiveddate{#1}}
\def\revised#1{\def\reviseddate{#1}}
\def\accepted#1{\def\accepteddate{#1}}

\def\coverauthors#1{\def\thecoverauthors{#1}}
\def\asciiauthors#1{\def\theasciiauthors{#1}}

\def\shortauthors#1{\def\theshortauthors{#1}}

\let\\\par\let\thevolumenumber\relax\let\thepapernumber\relax
\let\thevolumeyear\relax\let\thesamplenumber\relax\let\startpage\relax
\let\finishpage\relax\let\publishdate\relax\let\receiveddate\relax
\let\reviseddate\relax\let\accepteddate\relax\let\theasciititle\relax
\let\theasciiauthors\relax
\let\theasciiabstract\relax
\let\theasciiemail\relax\let\theshortauthors\relax\let\theshorttitle\relax
\let\thecoverauthors\relax

\long\def\maketitlep{   

\count0=\startpage

\gt\hfill      
\beginpicture
\setcoordinatesystem units <0.33truein, 0.33truein> point at 2.2 0.9
\setplotsymbol ({$\cal G$})
\plotsymbolspacing=9truept
\circulararc 315 degrees from 0 1 center at 0 0
\setplotsymbol ({$\cal T$})
\circulararc 315 degrees from 1 -1 center at 1 0
\endpicture
\break
{\small\ifx\thesamplenumber\relax 
Volume \else Sample
\fi\thevolumenumber\ (\thevolumeyear)
\startpage--\finishpage\nl
Published: \publishdate}
\vglue 0.5truein plus 0.4fil minus 0.1truein

{\parskip=0pt\leftskip 0pt plus 1fil\def\\{\par\smallskip}{\ifplaintex\large
\else\Large\fi\bf\thetitle}\par\medskip}   

\vglue 0pt plus 0.1fil 

{\parskip=0pt\leftskip 0pt plus 1fil\def\\{\par}{\sc\theauthors}
\par\medskip}

\vglue 0pt plus 0.1fil 

{\small\parskip=0pt\let\newline\\
{\leftskip 0pt plus 1fil\def\\{\par}{\sl\theaddress}\par}
\expandafter\ifx\theemail\relax    
\relax\else\vglue 5pt plus 0.02fil minus 2pt\def\\{\stdspace{\rm 
and}\stdspace} 
\cl{Email:\stdspace\tt\theemail}\fi
\ifx\theurl\relax                  
\relax\else\vglue 5pt plus 0.02fil minus 2pt\def\\{\stdspace{\rm 
and}\stdspace}
\cl{URL:\stdspace\tt\theurl}\fi\par}

\vglue 7pt plus 0.3fil minus 3pt

{\bf Abstract}
\vglue 5pt plus 0.1fil minus 2pt

\theabstract

\vglue 7pt plus 0.3fil minus 3pt

{\bf AMS Classification numbers}\quad Primary:\quad \theprimaryclass

Secondary:\quad \thesecondaryclass

\vglue 5pt plus 0.3fil minus 2pt

{\bf Keywords:}\quad \thekeywords

\vglue 10pt plus 0.5fil minus 5pt

{\small  Proposed: \theproposer\hfill Received: \receiveddate\nl
Seconded: \theseconders\hfill 
\ifx\reviseddate\relax                         
Accepted: \accepteddate                        
\else
Revised: \reviseddate                          
\fi}
\eject
}       

\let\maketitlepage\maketitlep
\let\maketitle\maketitlepage

\font\phead=cmsl9 scaled 950
\font=cmsl9 scaled 1050
\font\pnum=cmbx10 scaled 913
\font\lnum=cmbx10 
\font\pfoot=cmsl9 scaled 950
\font=cmsl9 scaled 1050
\ifplaintex
\headline{\vbox to 0pt{\vskip -4.5mm\line{\small\phead\ifnum
\count0=\startpage ISSN 1364-0380 (on line)
1465-3060 (printed) \hfill {\pnum\folio}\else\ifodd\count0\def\\{ }%
\ifx\theshorttitle\relax\thetitle\else\theshorttitle\fi\hfill{\pnum\folio}
\else\def\\{ and }{\pnum\folio}\hfill\ifx\theshortauthors\relax\theauthors
\else\theshortauthors\fi\fi\fi}\vss}}
\footline{\vbox to 0pt{\vglue 0mm\line{\small\pfoot\ifnum\count0=\startpage
\copyright\ \gtp\hfill\else
\gt, Volume \thevolumenumber\ (\thevolumeyear)\hfill\fi}\vss
}}
\else
\makeatletter
\def\@oddhead{{\small\lhead\ifnum\count0=\startpage ISSN 1364-0380 (on line)
1465-3060 (printed) \hfill {\lnum\number\count0}\else\ifodd\count0
\def\\{ }\ifx\theshorttitle\relax \thetitle \else\theshorttitle\fi\hfill
{\lnum\number\count0}\else\def\\{ and }{\lnum\number\count0}
\hfill\ifx\theshortauthors\relax 
\theauthors\else\theshortauthors\fi\fi\fi}}\def\@evenhead{@oddhead}
\def\@oddfoot{\small\lfoot\ifnum\count0=\startpage\copyright\ \gtp\hfill\else
\gt, Volume \thevolumenumber\ (\thevolumeyear)\hfill\fi}
\def\@evenfoot{@oddfoot}
\makeatother
\fi

\newwrite\gtoutfile
\long\gdef\makeheadfile{  
{\def\\{, }
\immediate\openout\gtoutfile head.xxx
\immediate\write\gtoutfile{To: math@arxiv.org}
\immediate\write\gtoutfile{Subject: put}
\immediate\write\gtoutfile{--text follows this line--}
\immediate\write\gtoutfile{Proxy-for: \ifx\theasciiauthors\relax
\theauthors\else\theasciiauthors\fi <\ifx\theasciiemail\relax\theemail\else\theasciiemail\fi>}
\immediate\write\gtoutfile{\noexpand\\}
\immediate\write\gtoutfile{Authors: \ifx\theasciiauthors\relax
\theauthors\else\theasciiauthors\fi}
{\def\\{ }\immediate\write\gtoutfile{Title: \ifx\theasciititle\relax
\thetitle\else\theasciititle\fi}}
\immediate\write\gtoutfile{Subj-class: GT}
\immediate\write\gtoutfile{MSC-class: \theprimaryclass\ifx\thesecondaryclass\relax\else, \thesecondaryclass\fi}
\immediate\write\gtoutfile{Journal-ref: Geom. Topol. \thevolumenumber
(\thevolumeyear) \startpage-\finishpage}
\immediate\write\gtoutfile{Comments: Published in Geometry and Topology at}
\immediate\write\gtoutfile{    http://www.maths.warwick.ac.uk/gt/GTVol\thevolumenumber/paper\thepapernumber.abs.html}
\immediate\write\gtoutfile{\noexpand\\}
\immediate\write\gtoutfile{}
\ifx\theasciiabstract\relax
\immediate\write\gtoutfile{\theabstract}\else
\immediate\write\gtoutfile{\theasciiabstract}\fi
\immediate\write\gtoutfile{}
\immediate\write\gtoutfile{\noexpand\\}
\immediate\write\gtoutfile{}
\immediate\write\gtoutfile{<uuencoded .tar.gz file here>}
\immediate\write\gtoutfile{}
\immediate\closeout\gtoutfile}}  

\def\maketitlepage{\maketitlep\makeheadfile}
\let\maketitle\maketitlepage

\volumenumber{4}\papernumber{7}\volumeyear{2000}
\pagenumbers{219}{242}
\proposed{Yasha Eliashberg}
\seconded{Robion Kirby, David Gabai}
\received{21 April 2000}
\accepted{1 September 2000}
\revised{20 July 2000}
\published{12 September 2000}

\newtheorem{thm}{Theorem}[section]

\newtheorem{prop}[thm]{Proposition}
\newtheorem{lemma}[thm]{Lemma}
\newtheorem{cor}[thm]{Corollary}

\newtheorem*{mainthm}{Theorem 6.1}
\theoremstyle{remark}
\newtheorem{defn}[thm]{Definition}

\newcommand{\R}{\mathbb{R}}
\newcommand{\Z}{\mathbb{Z}}

\newcommand{\bdry}{\partial}
\newcommand{\Cal}{\mathcal}
\newcommand{\sa}{\rightsquigarrow}

\newcommand{\be}{\begin{enumerate}}
\newcommand{\ee}{\end{enumerate}}
\nocolon

\begin{document}
\title{Tight contact structures and taut foliations}

\author{Ko Honda\\William H Kazez\\Gordana Mati\'c}
\asciiauthors{Ko Honda\\William H Kazez\\Gordana Matic}
\coverauthors{Ko Honda\\William H Kazez\\Gordana Mati\noexpand\'c}
\shortauthors{Honda, Kazez and Mati\'c}
\address{Mathematics Department, University of Georgia\\Athens, GA 30602, USA}
\email{honda@math.uga.edu, will@math.uga.edu, gordana@math.uga.edu}
\url{http://www.math.uga.edu/\char126 honda, http://www.math.uga.edu/\char126 will}



\begin{abstract}
We show the equivalence of several notions in the theory of taut
	foliations and the theory of tight contact structures. We
	prove equivalence, in certain cases, of existence of tight
	contact structures and taut foliations.
\end{abstract}
\maketitlepage

\section{Introduction}

The goal of this paper is to relate aspects of the
theory of taut foliations and the theory of tight contact structures. 
Codimension--1 foliations of 3--manifolds have a rich and beautiful history.
Highlights include the first examples on $S^3,$ due to Reeb, Haefliger's proof of the non-existence of
analytic foliations on $S^3$, and  Novikov's proof of the necessity of Reeb
components in foliations of $S^3$.  As a result of Gabai's work, the class of
foliations that have played the most important role in 3--dimensional 
topology, and especially in knot theory, are the taut foliations.
The theory of tight contact structures, on the other hand, has not yet reached a phase
where it can be applied effectively to the study of the topology of 3--manifolds.  It is still
concerned with basic questions about the structures themselves, such as existence and
classification on even some of the simplest manifolds, such as handlebodies. The
classification on $S^3$ and $B^3$ is due to Eliashberg in 1991, and the classification on
$T^3$ has been known only since 1995 \cite{K97, Gi94}. Only recently has the classification
been completed for lens spaces $L(p,q)$ \cite{Gi99, H1} and has the first
example of a manifold with no tight contact structure been produced \cite{EH99}
(the Poincar\'e
homology sphere $\Sigma(2,3,5)$ with one of its orientations).

Any relationship between these structures is not only interesting in its own right, but also
provides hope and an indication that contact structures will become a valuable
tool for studying 3--dimensional topology. Eliashberg and Thurston \cite{ET} bridged
the gap between foliation theory and contact topology.  Their seminal work opened the
door and enabled an exchange of ideas between two neighboring fields. They proved
that if a 3--manifold carries a taut foliation, then it also supports a tight contact
structure (in fact, one for each orientation of the ambient manfold $M$).
Although their method of perturbing a foliation into a contact structure is 3--dimensional,
their method of proving tightness is not 3--dimensional, and
instead uses the results from 4--dimensional
symplectic topology on symplectic fillings.
In \cite{HKM} we  reprove, and partially extend, their theorem using purely
3--dimensional techniques. The purpose of this paper is to prove a 
converse, in the
case of a 3--manifold with boundary, namely that if it supports a tight contact
structure, it supports a taut foliation. Note that we cannot hope to prove the
converse in the case of a general closed manifold,  since there are 
simple examples,
like $S^3$, which support tight contact structures but carry no taut 
foliations.

The techniques we use are based on a Haken decomposition theory, where the cutting
manifolds are {\it convex surfaces}.  In Section \ref{Convex surfaces} of this paper, we  briefly explain
the notion of a convex surface in a contact manifold as introduced by Giroux
\cite{Gi91}.  These appear to us to be the best kind of cutting surface for a
decomposition of a manifold with a contact structure.
In Section \ref{preparation}, we explain how to perturb a convex surface and a
(not necessarily Legendrian) curve $\gamma$ on it, so that $\gamma$ becomes Legendrian.
In Section \ref{Convex decompositions} we will explain how
to cut along convex surfaces with Legendrian boundary to perform a 
{\it convex splitting} on $M$. These will be used to cut the manifold eventually down to a union of
balls. Each ball supports a unique tight contact structure up to isotopy rel boundary, by a
fundamental theorem of Eliashberg \cite{E92}. The contact structure on $M$
is therefore encoded in the splitting surfaces $S$ together with characteristic foliation
on $S$. Moreover, the characteristic foliation on a convex surface $S$ is better encoded
by a collection of curves called the {\it dividing set } $\Gamma_S$.
Abstracting the idea of a 3--manifold $M$ with `curved' boundary $(\bdry M,\Gamma)$ ($\Gamma$
is a collection of curves), we define the notion of a {\it convex structure}. This notion
closely resembles the notion of {\it sutured manifolds} introduced by Gabai \cite{Ga} which we will
recall in Section \ref{sutured}. Gabai used {\it sutured manifold decompositions} to construct taut
foliations. We will show  that a convex Haken decomposition is, in a sense, a
generalization of a sutured manifold decomposition, and that the existence of a tight contact
structure on a manifold with given convex structure on the boundary implies the
existence of a taut foliation with the corresponding sutured manifold structure.  Our
main result, which incorporates important results of Gabai, Thurston and Eliashberg is:

\begin{mainthm} Let $(M,\gamma)$ be an irreducible sutured manifold
with annular sutures, and let $(M,\Gamma)$ be the associated convex
structure. The
following are equivalent.
\begin{enumerate}
\item[\rm(1)] $(M,\gamma)$ is taut.
\item[\rm(2)] $(M,\gamma)$ carries a taut foliation.
\item[\rm(3)] $(M,\Gamma)$ carries a universally tight contact structure.
\item[\rm(4)] $(M,\Gamma)$ carries a tight contact structure.
\end{enumerate}
\end{mainthm}

\section{Convex surfaces and convex structures}     \label{Convex surfaces}

Let $M$ be an oriented, compact 3--manifold (possibly with boundary).  A {\it
co-oriented positive contact structure} on $M$ is a nowhere 
integrable $2$--plane
field $\xi \subset T_*M$ such that there is a global 1--form $\alpha$ for which
$\alpha\wedge d\alpha = f \Omega$ with $f>0$ and $\Omega$ a volume form, and for which $\xi =
\rm ker\alpha$. $\alpha$ determines the orientation of $\xi$. A 
curve that is
everywhere tangent to the contact structure $\xi$ is called {\it 
Legendrian}. If
$\Sigma$ is an embedded surface, $\xi$ induces on it the  {\it characteristic
foliation} $\xi|_\Sigma$,which is defined to be the singular
foliation consisting of
the integral curves of $\xi \cap T_*\Sigma$ on $\Sigma$. Clearly these integral
curves are Legendrian.

A contact structure $\xi$ is said to be {\it overtwisted} if there exists a
disk $D$ which is everywhere tangent to $\xi$ along the boundary.  Such a disk $D$
is called an {\it overtwisted disk}.  A contact structure $\xi$ which is not overtwisted
is said to be {\it tight}.  Eliashberg \cite{E89} showed that, for closed 3--manifolds,
the set of overtwisted contact 2--plane fields is weak homotopy equivalent to
the set of contact 2--plane fields (without any integrability conditions).  Hence, the study of
overtwisted contact structures is largely homotopy-theoretic (of course there is the problem of determining
whether a contact structure is tight or overtwisted). Tight contact structures are less
ubiquitous, and tend to reflect the topology of the 3--manifold in ways which are not very
well-understood.

We say a vector field $v$ on a contact manifold $(M,\xi)$ is  a {\it contact vector field} if
its flow preserves $\xi$. An oriented properly embedded surface
$\Sigma$ in $(M,\xi)$ is called {\it convex} if there is a contact vector
field $v$ transverse to $\Sigma.$  The {\it dividing set} $\Gamma_\Sigma$
of a convex surface $\Sigma$ with respect to a transverse contact vector
field $v$ is the set of points $x$ for which $v(x)\in \xi(x)$.
The following is a fundamental theorem of Giroux \cite{Gi91}.

\begin{thm} [Giroux \cite{Gi91}]\label{t:fundamental} The dividing set
$\Gamma_\Sigma$ is a union of smooth curves which are transverse to the
characteristic foliation $\xi|_\Sigma$. Moreover, the isotopy type of 
$\Gamma_\Sigma$
is independent of the choice of $v$.
\end{thm}

The isotopy class of
$\Gamma_\Sigma$ is clearly preserved under an isotopy of $\Sigma$
through a family of convex surfaces. Conversely, if ${\cal
F}$ is a singular foliation on $\Sigma$, then a disjoint union of properly
embedded curves $\Gamma$ is said to {\it divide} ${\cal F}$ if there exists
an $I$--invariant contact structure $\xi$ on $\Sigma\times I$ such that
${\cal F}=\xi|_{\Sigma\times \{0\}}$ and $\Gamma$ is the dividing set for
$\Sigma\times \{0\}$.

Denote the
number of connected components of $\Gamma_\Sigma$ by $\#\Gamma_\Sigma$. The
complement of the dividing set is the union of two subsets
$\Sigma\backslash\Gamma_\Sigma = R_+-R_-$. Here $R_+$ is the subsurface
where the orientations of $v$ and the normal orientation of $\xi$ coincide,
and $R_-$ is the subsurface where they are opposite. If $\Sigma$ is a
surface with boundary, in this paper we also require that the
boundary be a Legendrian curve for $\Sigma$ to be called {\it convex}.

\begin{thm}[Giroux's Flexibility Theorem \cite{Gi91}]\label{t:flexibility}
Let $\Sigma$ be a convex surface in a contact $3$--manifold $(M, \xi)$ with
characteristic foliation $\xi|_\Sigma$, contact vector field $v$, and
dividing set $\Gamma$. If $\cal{F}$ is another singular foliation on
$\Sigma$ divided by $\Gamma$, then there is an isotopy $\phi_t\co \Sigma\rightarrow M$, $t\in[0,1]$,
such that $\phi_0(\Sigma)=\Sigma,$
$\xi|_{\phi_1(\Sigma)}=\cal{F}$, the isotopy is fixed on
$\Gamma$, and $\phi_t(\Sigma)$ is transverse to $v$ for all $t$.
\end{thm}

Such an isotopy is said to be an {\it admissible  isotopy} of a convex surface $\Sigma$
with respect to a contact vector field $v\pitchfork \Sigma$.   If the contact vector field
$v$ is omitted, it is implied that the isotopy is admissible with respect to some $v$.

Giroux also finds conditions under which a convex surface has a tight
$I$--invariant contact neighborhood.

\begin{thm}[Giroux]\label{t:convextight} If $\Sigma\not = S^2$ is a
convex surface in a contact manifold $(M,\xi)$, then $\Sigma$ has a tight
neighborhood if and only if no component of $\Gamma_\Sigma$ is null-homotopic
in $\Sigma$. If $\Sigma=S^2$, $\Sigma$ has a tight neighborhood if and only if
$\#\Gamma_\Sigma=1$.
\end{thm}

We say that a contact structure on a manifold $M$ with boundary
$\partial M$ is a {\it contact structure with convex boundary} if 
there is a contact
vector field $v$ on $M$ transverse to $\partial M$. The following 
definition records the information about a contact structure near
its convex boundary, but forgets the structure in the interior.

\begin{defn}  A {\it convex structure} is a quadruple $(M, \Gamma,
R_-(\Gamma), R_+(\Gamma))$ where $M$ is a compact oriented 3--manifold
with nonempty boundary,
$\Gamma$ is a disjoint union of simple closed curves contained in
$\partial M$ nonempty on each component of $\partial M$, and $\partial M =
R_+(\Gamma) \cup R_-(\Gamma)$, $R_+(\Gamma) \cap R_-(\Gamma) = \Gamma$.
Moreover  $R_+(\Gamma), R_-(\Gamma)$ and $\Gamma$ are oriented so that the
orientation of $R_+(\Gamma)$ agrees with the orientation induced on $\partial
M$ by the orientation of $M$, and the orientation on $R_-(\Gamma)$ is the
opposite one. $\Gamma$ is oriented in such a way that if $\alpha \subset
\partial M$ is an oriented arc with $\partial \alpha \subset R_+(\Gamma) \cup
R_-(\Gamma)$ that intersects $\Gamma$ transversely in one point and if $\Gamma
\cdot \alpha = 1$ then $\alpha$ must start in $R_-(\Gamma)$ and end in
$R_+(\Gamma)$. \end{defn}

 A contact structure on $M$ with convex boundary and a
choice of a contact vector field $v$ such that $v$ is an oriented normal to
$\partial M$ induces a convex structure on
$M$. $\Gamma$ is defined to be the dividing set of $\bdry M$ with respect to $v$, and $R_+(\Gamma)$
and $R_-(\Gamma)$ are the regions of $\partial M$ where the oriented normal
vector $n_{\xi}$ to the contact planes  and $v$ satisfy $n_{\xi} \cdot v >0$
and $n_{\xi} \cdot v <0$ respectively.

\begin{defn} A convex structure $(M,\Gamma, R_-(\Gamma), R_+(\Gamma))$
{\it carries a tight contact structure} if there is a tight contact
structure on $M$, and a contact vector field $v$ such that
$v$ is an oriented normal for $\partial M$
and both $\Gamma, R_-(\Gamma)$ and $R_+(\Gamma)$ are defined by $v$ as
above.
\end{defn}

Note that if we change the orientation of the contact plane field $\xi$,
$R_-(\Gamma)$ and $R_+(\Gamma)$  will
switch.

\section{Legendrian curves on convex surfaces}         \label{preparation}

A Legendrian curve $C$ and the oriented normal to $\xi$ determine a framing
along $C$. If $Fr$ is another framing we define the {\it twisting number}
$t(C,Fr)$ as the relative framing
between the one determined by the oriented normal to $\xi$ and  $Fr$.  If $C$
lies on a surface $\Sigma$, $t(C,\Sigma)$ is defined to be the twisting number
with respect to the framing defined on $C$ by  $\Sigma$. Observe that
if if $C$ is a
Legendrian curve on a convex surface $\Sigma$, then its twisting number
$t(C,\Sigma)$ is equal $\frac{1}{2}\#(C\cap \Gamma_\Sigma)$, where $\#(C\cap
\Gamma_\Sigma)$ denotes the geometric intersection number. In fact it is easy
to show the following.

\begin{prop} Let $C$ be a Legendrian curve on a convex surface $\Sigma$ with
$t(C,\Sigma)=-n$. Then, after a small perturbation of $\Sigma$,
there are local coordinates $(x,y,z)$ so that a
neighborhood of  $C$ in $M$
is isomorphic to the neighborhood  $N=\{(x,y,z)|x^2 + y^2  \le
\varepsilon\}
$ in $\R^2 \times (\R/\Z)$,
where the set $x=0$ corresponds to $\Sigma$, $C$ is given by $x=y=0$,
and the contact structure is
determined by the 1--form $\alpha = \sin(2\pi n z) dx + \cos(2\pi n z) dy $.
If the contact vector field determining the dividing set $\Gamma_\Sigma$ is
$v=\frac{\partial}{\partial x}$, the dividing set is $\Gamma_\Sigma =
\{(0,y,\frac{k}{2n})| 0 \le k \le 2n\}$.
\end{prop}

It is a standard fact that any curve in a contact manifold has in its
isotopy class
a nearby Legendrian curve. However, even more is true: this can be
achieved even when we require the curve to
lie on a convex surface isotopic to a fixed one and with the same dividing set.
Let us call a union of closed curves
$C$ on a convex surface $\Sigma$  {\it nonisolating} if
(1) $C$ is transverse to $\Gamma_\Sigma$,  and (2) every component of
$\Sigma\backslash (\Gamma_\Sigma \cup C)$ has a boundary component which
intersects $\Gamma_\Sigma$. Clearly this will be satisfied if every
component of $C$ intersects $\Gamma_\Sigma$.

\begin{thm}[Legendrian Realization Principle \cite{H1}] \label{t:lerp}
Let $C$ be a nonisolating collection of closed curves on a convex
surface $\Sigma$.  Then there exists an admissible isotopy $\phi_t$,
$t\in[0,1]$, so
that
\be
\item $\phi_0=id$,
\item $\phi_t(\Sigma)$ are all convex,
\item $\phi_1(\Gamma_\Sigma)=\Gamma_{\phi_1(\Sigma)}$,
\item $\phi_1(C)$ is Legendrian.
\ee
\end{thm}

It follows that a nonisolating collection $C$ can be realized by a
Legendrian collection $C'$ with the same number of geometric intersections with
$\Gamma_\Sigma$.  A special case of this theorem, observed by Kanda, is
the following:

\begin{cor}[Kanda] If $C$ is a closed curve in $\Sigma$ such that $C\pitchfork \Gamma_\Sigma$
and  $C \cap
\Gamma_\Sigma \ne \emptyset$,
then $C$ can be realized as a Legendrian curve (in the sense of
Theorem~\ref{t:lerp}).
\end{cor}

Giroux \cite{Gi91} proved that a closed oriented embedded surface can 
be deformed
through a $C^\infty$--small isotopy to a convex surface.  The following relative
version is proven in Honda \cite{H1}.

\begin{thm}[Existence of Convex Surfaces] \label{existence}  Let $T 
\subset M$ be a
compact, oriented, properly embedded surface with Legendrian boundary 
such that $t(C,
T)$ $\le 0$ for all components $C$ of $\partial T$.  There exists a 
$C^0$--small isotopy
of $T$, which is the identity on $\partial T$, that takes $T$ to a 
convex surface.
The isotopy may be chosen to be $C^\infty$ outside of a small neighborhood of
$\partial T$.
\end{thm}

\section{Convex decompositions}       \label{Convex decompositions}

A 3--manifold $M$ is {\it irreducible} if every embedded 2--sphere $S^2$ bounds a
3--ball $B^3$.     A properly embedded surface $\Sigma\subset M$ is {\it incompressible}
if it contains no {\it compressing disk}, ie, an embedded disk $D\subset M$ with
$D\cap \Sigma=\bdry D$ which is homotopically nontrivial in $\Sigma$.
 A {\it Haken decomposition} of a 3--manifold $M$ is a sequence
\begin{equation}  \label{Haken}
M=M_0\stackrel{S_1}{\sa} M_1\stackrel{S_2}\sa \cdots \stackrel{S_n}\sa M_n,
\end{equation}
where $ S_{i+1} $ is an incompressible surface in $M_i$, $M_{i+1} = M_i
\backslash S_{i+1} $, and $M_n$ is a disjoint union of balls.
{\it Haken manifolds} are 3--manifolds which admit Haken decompositions.
Therefore, inductive arguments can often be applied to Haken manifolds.  An
irreducible manifold with non-empty boundary always has a Haken decomposition
\cite{J}. The idea we are pursuing in this paper is that when $M$
has a contact
structure, and we choose the splitting surfaces to be convex, the information about
the contact structure on $M$ can be recovered from the contact 
structure on the cut-up manifold $M \backslash S$ and the information contained in the
dividing set on the splitting surface $S$.   In this section we will describe how to
perform {\it convex splittings} in the contact category.

When $(M, \xi)$ is a contact structure with convex boundary, we can choose
a Haken decomposition of $(M, \partial M)$ to be, at each step, performed along
incompressible surfaces  with boundary $(S,\partial S)$ properly embedded in
$(M,\partial M)$. At each step of the decomposition, we will do the
same three things: perturb the cutting surface $(S,\partial S)$ to a convex
surface with Legendrian boundary, cut $(M, \xi)$ along $S$ to obtain a
manifold with corners $M \backslash S$ which inherits the restriction
$\xi|_{M\backslash S}$ of $\xi$, and finally round corners to obtain a smooth manifold
and a contact structure with convex boundary on it.

\begin{figure}[ht!]
	{\epsfysize=2in\centerline{\epsfbox{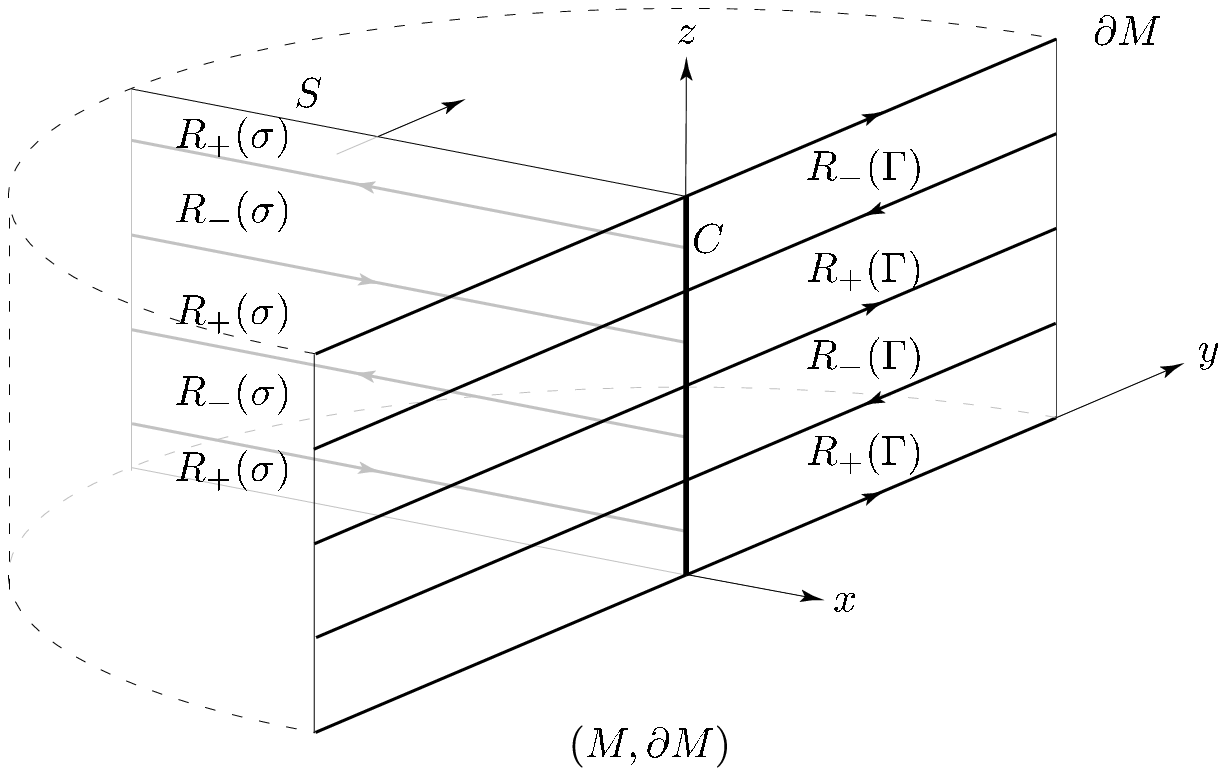}}}
	\caption{}
	\label{f:local}
\end{figure}

We first need to perturb $\partial S$. We isotop each component $C$ of
$\bdry S \subset \partial M$ so that the geometric intersection $\#(C\cap
\Gamma_ {\partial M} )$
 is minimized, provided this number is $ \geq 2$.  If the minimum
geometric intersection is $0$, we can choose $C$ so $\#(C\cap \Gamma_{\partial M} )
=2$, since every component of $\bdry M$ nontrivially intersects $\Gamma_{\bdry M}$.
We artificially force the extra intersections
because cutting along Legendrian curves with twisting number $0$ is
not as easy to control. Now we can
use the Legendrian Realization Principle (Theorem~\ref{t:lerp}) to make
$\bdry S$ Legendrian.  Once we have prepared $\partial S$ as above, we
perturb the surface $S$ so
that near the boundary it is convex and the local picture is as
in Figure~\ref{f:local}.

If $C$ intersects the dividing set $\Gamma_{\partial M}$ geometrically $2n$
times, there is a neighborhood of $C$ in $M$ and local coordinates $(x,
y, z)$ on it isomorphic to $N = \{(x, y,z)|x^2 + y^2 < \varepsilon, x \le
0\}$ in $\R^2 \times (\R/\Z)$ where the set
$A = \{(x,y,z) \in N |x = 0\}$
corresponds to an annular neighborhood of $C$ in $\partial M$ and $B
= \{(x,y,z) \in N |y = 0\}$ to an annular neighborhood of $C$ in $S$, and the
1--form $\alpha = \sin(2\pi n z) dx + \cos(2\pi n z) dy$ determines the contact
structure.  If we choose the contact vector fields for $\partial M$ and
$S$ in these coordinates to be respectively $v_{\partial
M}=\frac{\partial}{\partial x}$ and $v_S = \frac{\partial}{\partial y}$ it is
easy to calculate that the dividing sets are
$\Gamma_{\partial M} = \{0, y, \frac{k}{2n})| 0 \le k < 2n\}$ and
$\Gamma_S = \{(x, 0), \frac{1 + 2k}{4n}| 0 \le k < 2n\}$.

If $(M,\Gamma, R_+,R_-)$ is the convex
structure associated to a contact structure $\xi$ with convex boundary, and if
$S$ is a convex surface with Legendrian boundary
properly embedded in $M$ and transverse to $\Gamma$, then the convex
vector field  $v_S$ given by  $ \frac{\partial}{\partial y}$ in the 
local coordinates
discussed above can be extended to a convex vector field on $S$, which will
determine a dividing set $\sigma$ on $S$
as well as  subsets
$R_-(\sigma)$ and
$R_+(\sigma)$, defined as in the case of a closed surface.

The next definition abstracts the properties of a properly embedded convex
surface with Legendrian boundary in a contact manifold with convex boundary.

\begin{defn}  A {\em surface with divides} $(S, \sigma,
R_+(\sigma),R_-(\sigma) )$ is a
compact oriented surface $S$, possibly with boundary, together with a disjoint
collection of properly embedded arcs and simple closed curves $\sigma$ and a
decomposition into two subsurfaces $S=R_+(\sigma)\cup R_-(\sigma)$,
$R_+(\sigma)
\cap R_-(\sigma)= \sigma$.  The orientation on $R_+(\sigma)$ is the orientation
induced from $S$ while $R_-(\sigma)$ has the opposite orientation.
The components
of $\sigma$ are oriented so that if $\alpha \subset S$ is an oriented arc which
intersects $\sigma$ transversely in one point and $\sigma \cdot \alpha = 1$
then $\alpha$ starts in $R_-(\sigma)$ and ends in $R_+(\sigma)$.
\end{defn}

Dividing curves on convex surfaces in tight contact manifolds satisfy special
properties, as we saw in Theorem ~\ref{t:convextight}. For a convex surface
with Legendrian boundary we have the following generalization:

\begin{prop} Let $(M,\xi)$ be a tight contact manifold with convex boundary,
and  let $\sigma$ be the dividing set of a convex surface $S$
with Legendrian boundary $\partial S$ transverse to the dividing set
$\Gamma_{\partial M}$, such that every component of $\partial S$
intersects $\Gamma_{\partial M}$. Then  $\sigma$ satisfies the following:
\begin{enumerate}

\item On each component of $\partial S$ the points of $\sigma \cap
\partial S$ alternate with the points of $\Gamma \cap \partial S$.

\item The orientation on each arc of $\sigma$ is from $R_-(\Gamma)$
to $R_+(\Gamma)$.

\item No  closed curve in $\sigma$ bounds a disk in $S$.

\end{enumerate}
\end{prop}

\proof

Parts $1$ and $2$ follow from the local coordinates picture discussed above
and part $3$ from Theorem~\ref{t:convextight}.    \endproof


When we split $(M, \partial M)$ along $(S, \partial S)$ we obtain a
manifold with corners $M \backslash S$.  To smooth the corners we
use the following ``corner-rounding'' procedure. Each of the halves of
$N$,
$$N_- = \{(x, y,z) \in N| y \le 0\}$$
and
$$N_+ = \{(x, y,z) \in N|y \ge 0\}$$
is replaced by the corresponding
$$N^r_- = \Bigl\{(x, y,z) \in N_-| x \le \frac{-\varepsilon}{2} 
\rm{\ \  or  \ \
}y \le -\frac{\varepsilon}{2} + \sqrt{\frac{\varepsilon^2}{4} - (x
+\frac{\varepsilon}{2})^2}\  \Bigr\}$$
and
$$N^r_+ = \Bigl\{(x, y,z)  \in N_+|x \le \frac{-\varepsilon}{2} \rm{\ \ or
\ \ }y \le +\frac{\varepsilon}{2} + \sqrt{\frac{\varepsilon^2}{4} - (x
+\frac{\varepsilon}{2})^2}\  \Bigr\}.$$

\begin{figure}[ht!]

{\epsfysize=2in\centerline{\epsfbox{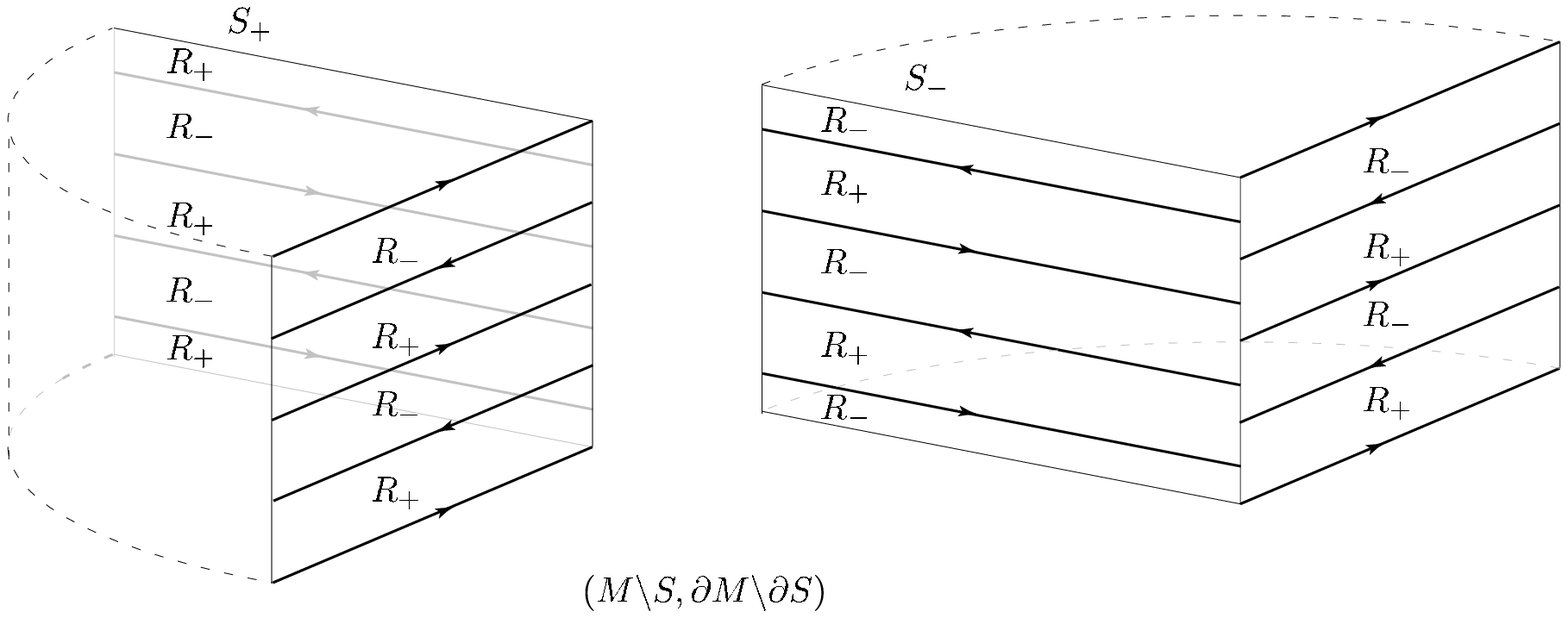}}}
	\caption{}
	\label{SPLITTING}
\end{figure}

A quick look at the form $\alpha = \sin(2\pi n z) dx + \cos(2\sin z) dy$
determining $\xi$ and the normal vectors of the boundaries show, even without
calculation, that the dividing set on the rounded boundary will be as in
Figure 3. Clearly, $\mbox{ker} \alpha = \mbox{span}\{\frac{\partial}{\partial
z},  \cos(2 \pi nz)\frac{\partial}{\partial x} - \sin(2 \pi
nz)\frac{\partial}{\partial y}\}$, and the contact vector fields all lie in the
$(x,y)$--plane. It is an easy calculation to see that when the contact vector
rotates counterclockwise in the $(x,y)$--plane, the $z$--coordinate of the
dividing set decreases.

\begin{figure}[ht!]

{\epsfysize=2in\centerline{\epsfbox{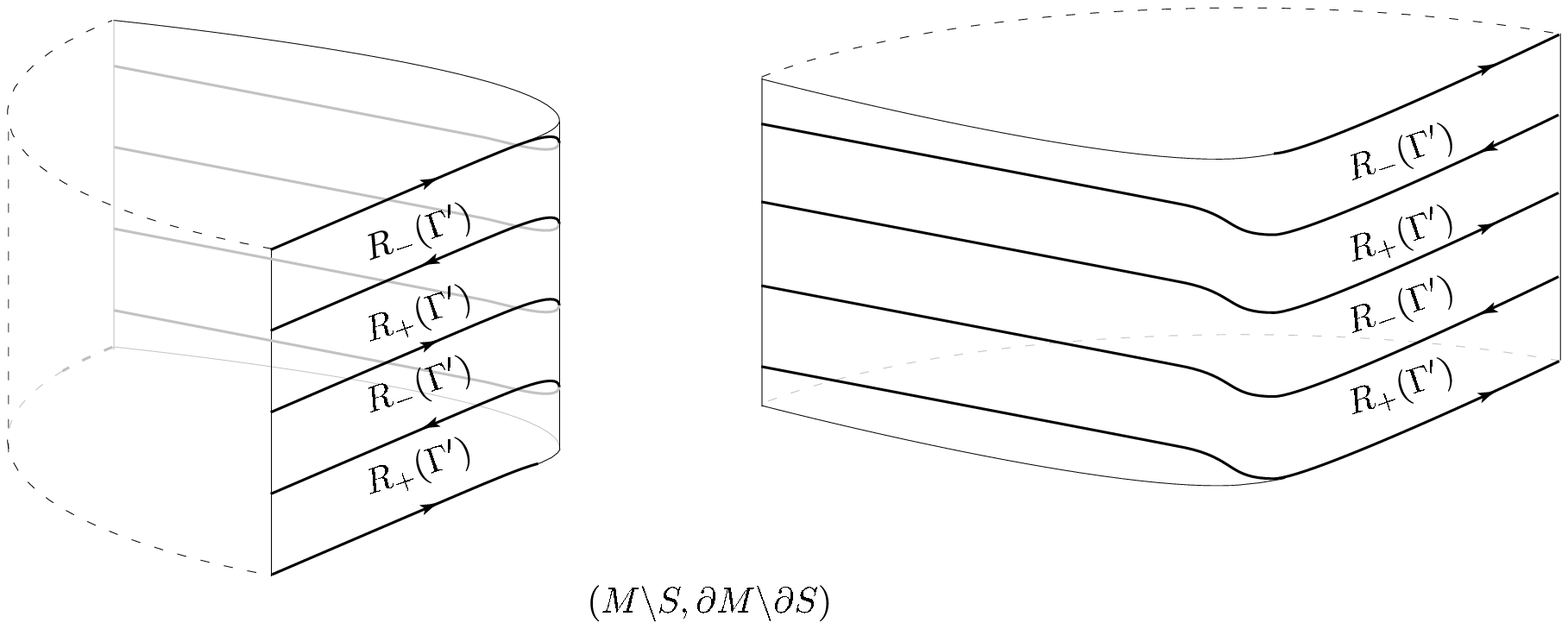}}}
	\caption{}
	\label{ROUNDING CORNERS}
\end{figure}

We introduce the notion of a {\it convex splitting} to formalize the
proces of obtaining the convex structure on the manifold with boundary
$(M
\backslash S, \partial M \backslash \partial S)$ by cutting $(M,\partial M)$
along the properly embedded convex surface with Legendrian boundary $S$,
rounding the corners and looking at the new dividing set.

\begin{defn} Let $(S,\sigma)$ be a surface with divides that is
properly  embedded
in a convex structure $(M,\Gamma)$ so that $S$ and $\sigma$ are  both
transverse
to $\Gamma$, and so that they satisfy properties 1--\,3 listed  above. We say
that $(S,
\sigma)$ defines a {\it convex splitting} $(M, \Gamma) \stackrel{(S,
\sigma)}{\rightsquigarrow} (M', \Gamma')$.   $M'$ is $M$ split along $S$ and is
denoted
$M' = M\backslash S$.  $\partial M'$ contains two  disjoint copies of $S$ which
are denoted $S_+$ and $S_-$.  $S_+$ are the components such that the  outward
orientation it inherits from $M'$ agrees with the original orientation on $S$.
Given a subset $X \subset S$ denote by $X_+$ the corresponding subset of $S_+$,
and  similarly for $X_-$.  Thus $\sigma_+, (R_+(\sigma))_+,
(R_-(\sigma))_+$ are
all  subsets of  $S_+$. Define \begin{eqnarray*}
 R_+(\Gamma') &= & (R_+ (\Gamma)\backslash \partial S) \cup
(R_+(\sigma))_+ \cup
(R_-(\sigma))_-\\ R_-(\Gamma') &= & (R_-(\Gamma)\backslash \partial S) \cup
(R_-(\sigma))_+ \cup (R_+(\sigma))_-\\ \Gamma' & = & R_+(\Gamma') \cap
R_-(\Gamma'). \end{eqnarray*} Finally, smooth all corners so that
$\partial M'$ is
a smooth subset of $M'$ and $\Gamma'$ is a smooth subset of $\partial M'$.
\end{defn}

\begin{figure}[ht!]

{\epsfysize=3.5in\centerline{\epsfbox{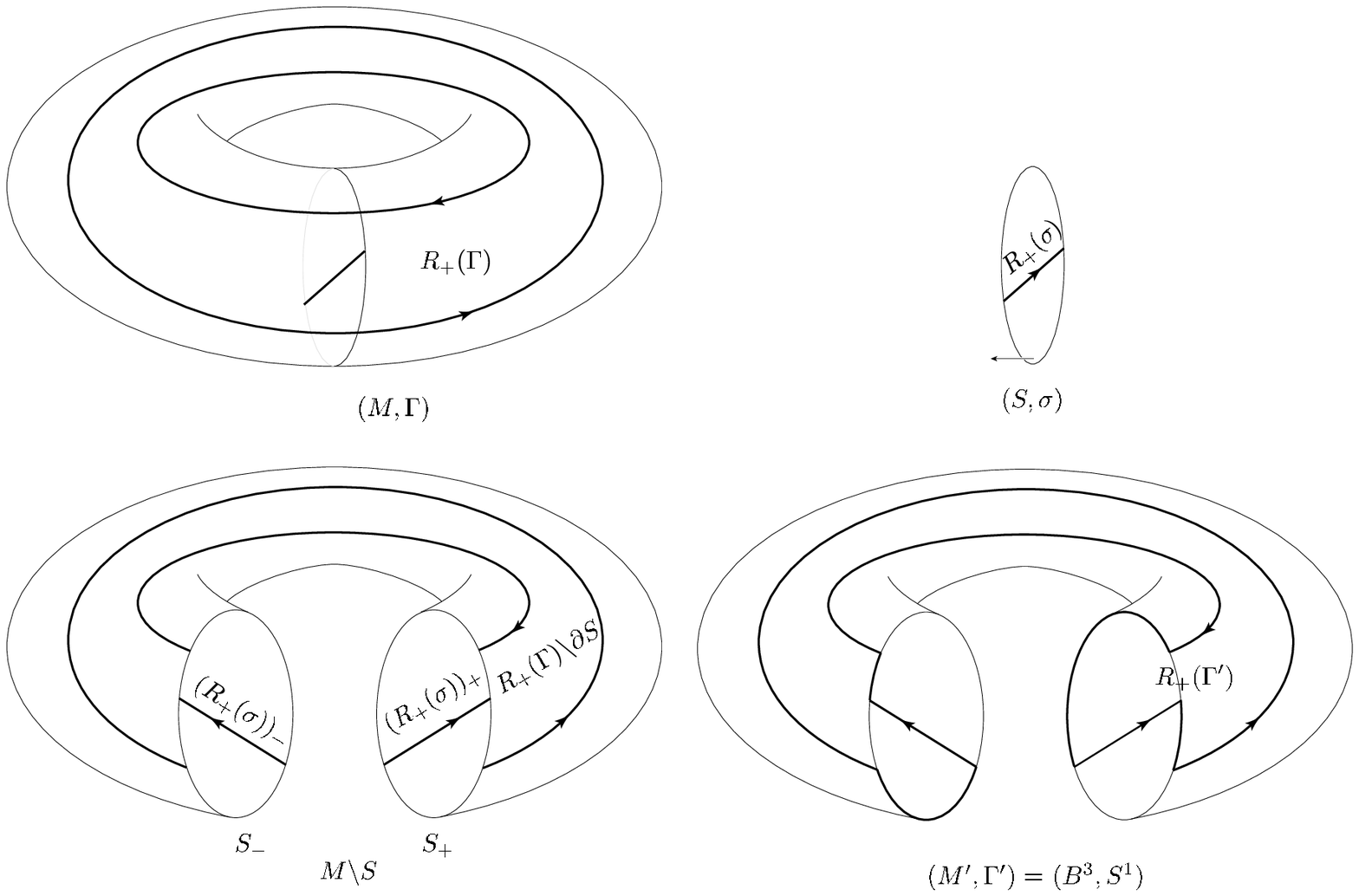}}}
	\caption{}
	\label{convexsplitting}
\end{figure}

If we perform a Haken decomposition of a tight contact manifold with convex
boundary along embedded convex surfaces with Legendrian boundary,
rounding corners
at each step along the way, we obtain in the end a disjoint union of spheres
with tight contact structures on them.   The following facts now come into play:

\begin{prop}  \label{pikachu} Let $\xi$ be a tight contact structure on $B^3$ with convex boundary.
Then $\#\Gamma_{\bdry B^3}=1$.
\end{prop}

This is just Theorem \ref{t:convextight} restated.

\begin{thm} [Eliashberg \cite{E92}]   Let $\xi$ be a contact structure on a neighborhood of
$\bdry B^3$ for which $\bdry B^3$ is convex and $\#\Gamma_{\bdry B^3}=1$.
Then there exists a unique extension of $\xi$ to a tight contact structure on
$B^3$, up to an isotopy which fixes the boundary.
\end{thm}

The decomposition of
tight contact manifolds motivates the following definition of
decomposability of
convex structures.

\begin{defn} A convex structure $(M, \Gamma)$ is {\it decomposable} if there
exists a sequence of convex splittings
$$(M, \Gamma) \stackrel{(S_1,\sigma_1)}{\rightsquigarrow} (M_1, \Gamma_1)
\rightsquigarrow \cdots \stackrel{(S_n,\sigma_n)}\rightsquigarrow
(M_n, \Gamma_n)$$
such that $(M_n, \Gamma_n)$ is a disjoint union of $(B^3, S^1)$'s.
\end{defn}

We then have the following:

\begin{thm}\label{A} If $(M, \Gamma)$ carries a tight contact
structure, then it is decomposable.
\end{thm}

\proof Let $\xi$ be a tight contact structure on $M$ which is adapted to $\Gamma$.
Consider the Haken decomposition
$$M=M_0\stackrel{S_1}{\sa} M_1\stackrel{S_2}\sa \cdots \stackrel{S_n}\sa M_n.$$
Let $\Gamma_0=\Gamma$.
Assume we have already performed convex splittings along convex surfaces with Legendrian
boundary, so that we have $(M_i,\Gamma_i)$.  In order to split along $S_{i+1}$ in a convex manner,
make $\bdry S_{i+1}$ Legendrian using the Legendrian Realization Principle, perturb $S_{i+1}$ so
it is convex with Legendrian boundary, form $M_i\backslash S_{i+1}$, and round the corners.
This yields $(M_{i+1},\Gamma_{i+1})$.  Since $M$ is Haken, we eventually find that
$M_n=\cup B^3$.  Proposition \ref{pikachu} implies that for each $B^3$ we have $\#\Gamma_{\bdry B^3}=1$.
\qed

\begin{cor}\label{A1} If $(M, \Gamma)$ carries a tight contact structure, then
$\chi(R_+(\Gamma))= \chi(R_-(\Gamma))$.
\end{cor}

\proof If $(M, \Gamma) \stackrel{(S,
\sigma)}{\rightsquigarrow} (M', \Gamma')$, then a computation
shows that $\chi(R_\pm(\Gamma')) = \chi(R_\pm(\Gamma)) + \chi(S)$.
The result follows by
induction on the length of the decomposition sequence for $(M, \Gamma)$. \qed

\section{Sutured vs convex decompositions}    \label{sutured}

We now recall basic definitions from Gabai's theory of sutured manifolds
\cite{Ga}. It will be immediately obvious that they resemble the
definitions just
made. The point of this paper is to exploit the equivalence of basic notions in
these theories.

\begin{defn}  A {\it sutured manifold} $(M, \gamma)$ is a compact oriented
3--manifold $M$ together with a set $\gamma \subset \partial M$ of pairwise
disjoint annuli $A(\gamma)$ and tori $T(\gamma)$.  $R(\gamma)$
denotes $\partial M\backslash int(\gamma)$.  Each component of $R(\gamma)$ is oriented. $R_+(\gamma)$ is
defined to be those components of $R(\gamma)$ whose normal vectors point out of
$M$ and $R_-(\gamma)$ is defined to be $R(\gamma) \backslash R_+(\gamma)$.
Each component of $A(\gamma)$ contains a {\it suture}, ie, a homologically
nontrivial oriented simple closed curve.  The set of sutures is denoted $s(\gamma)$.  The
orientation on $R_+(\gamma), R_-(\gamma)$ and $s(\gamma)$ are related as follows.
If $\alpha \subset \partial M$ is an oriented arc with $\partial \alpha \subset R(\gamma)$
that intersects $s(\gamma)$ transversely in one point and if $s(\gamma) \cdot
\alpha = 1$, then $\alpha$ must start in $R_-(\gamma)$ and end in
$R_+(\gamma)$. \end{defn}

\begin{defn}  A {\it sutured manifold with annular sutures}
is a sutured manifold
$(M, \gamma)$ which satisfies the following:
\be
\item Every component of $M$ has nonempty boundary.
\item Every component of $\bdry M$ contains a suture.
\item Every component of $\gamma$ is an annulus.
\ee
\end{defn}

Note that a sutured manifold $(M, \gamma)$ with annular sutures
determines, and is
determined by, the {\it associated convex structure} $(M, \Gamma)$
where $\Gamma = s(\gamma)$.

The definition of a {\it sutured manifold splitting} $(M,\gamma)
\stackrel{S}{\rightsquigarrow} (M',\gamma')$ is quite similar to the
definition of a convex splitting. However, unlike convex splittings,
we do not have dividing curves to prescribe on the splitting surface $S$.

Assume $S$ is a properly
embedded, oriented surface in $M$ such that:

\begin{enumerate}
\item $\bdry S \pitchfork \gamma$.

\item If $S$ intersects an annular suture $A$ in arcs, then no such
arc separates
$A$.

\item If $S$ intersects an annular suture $A$ in circles, then each such circle,
with orientation induced from $S$, is homologous in $A$ to the oriented core
$s(\gamma) \cap A$.

\item If $S$ intersects a toroidal suture $T$ in circles, then no such circle is
null-homologous in $T$, and any two such circles, with orientations
induced from$S$, are homologous in $T$.

\item No component of $S$ is a disk $D$ with $\partial D \subset R(\gamma)$.

\item No component of $\partial S$ bounds a disk in $R(\gamma)$.

\end{enumerate}

Let $M'=M\backslash S$ and let $S_+$ and $S_-$ be the copies of $S$
contained in $M'$
where the orientation induced by $S$ points, respectively, out of and
into $M'$.  As a
first approximation, let $R'_{\pm}(\gamma')$ be $(R_{\pm}(\gamma)\backslash
S)\cup S_{\pm}$.  $\gamma'$ is supposed to separate $R'_+(\gamma')$
and $R'_-(\gamma')$ so
define it to be the union of $\gamma\backslash S$ and
$R'_+(\gamma')\cap R'_-(\gamma')$.
Since $\gamma'$ is supposed to be a union of annuli and tori, the
actual definition of
$\gamma'$ is a union of $\gamma\backslash S$ and a regular neighborhood of
$R'_+(\gamma')\cap R'_-(\gamma')$ and then $R'_{\pm}(\gamma')$ are shrunk by a
corresponding amount.

\begin{defn}   A transversely oriented codimension--1 foliation $\cal F$ is
{\it carried by} $(M, \gamma)$ if $\Cal F$ is transverse to $\gamma$ and tangent to
$R(\gamma)$ with the normal direction pointing outward along $R_+(\gamma)$ and
inward along $R_-(\gamma)$, and $\Cal F|\gamma$ has no Reeb components. $\Cal F$ is
{\it taut} if each leaf of $\Cal F$ intersects some closed curve or properly embedded arc  connecting
from $R_-(\gamma)$ to $R_+(\gamma)$
that is transverse to $\Cal F$.
\end{defn}

Let $S$ be a compact oriented surface with components $S_1, \dots,
S_n$. The {\it
Thurston norm of $S$} is defined to be
$$
x(S) = \sum_{i\rm{\ such\ that\ } \chi(S_i) < 0} |\chi(S_i)|.
$$
Thus components with positive Euler characteristic, namely disks and
spheres, do
not contribute to the Thurston norm.

\begin{defn} A sutured manifold $(M, \gamma)$ is {\it taut} if

\begin{enumerate}
\item  $M$ is irreducible.

\item $R(\gamma)$ is {\it norm-minimizing} in $H_2(M, \gamma)$, that is if
$S$ is an embedded surface in $M$ with $[S] = [R(\gamma)] \in H_2(M,
\gamma)$ then $x(R(\gamma)) \le x(S)$.

\item $R(\gamma)$ is incompressible in $M$.
\end{enumerate}
\end{defn}

Except in a few cases, 2 implies 3.  The reason is that compressions
are norm-decreasing
unless the surface being compressed is an annulus.  Thus 3 is meant
to exclude the
case that $M = B^3$ and $s(\gamma)$ consists of more than one
component or that $M
= D^2 \times S^1$ and $s(\gamma)$ is compressible.

This definition of tautness of the sutured manifold is made because of the
following theorem which is due to Gabai \cite{Ga} and Thurston \cite{Th}.

\begin{thm} A sutured manifold $(M, \gamma)$ is taut if and only if
it carries a
transversely oriented, taut, codimension--1 foliation $\Cal F$.
\end{thm}

The following correspondence shows that a sutured manifold splitting is a
special case of the convex splitting:
\be
\item The cores of annular components of $\gamma$ can be viewed as  dividing
curves. If T is a toroidal component of $\gamma$ then just before cutting
along a surface $S$ which intersects $T$ we substitute  $T$  by $T$ with a
pair of parallel homotopically nontrivial dividing curves, each of  which has
algebraic intersection $ 1 $ with each component of $S \cap T$.

\item A component $\Sigma\subset \bdry M$ may not have a suture at
all, whereas a
dividing set must not be empty.  We remedy this by placing a pair of parallel
homotopically nontrivial dividing curves on $\Sigma$ before cutting.

\item Let  $S$ be a cutting surface -- realize the boundary as a
Legendrian curve
with twisting number $\leq -2$ -- and choose $\Gamma_S$ so that every dividing
curve is an arc which is $\bdry$--compressible.

\item When $M$ is cut along $S$ and rounded, all the dividing curves, except
perhaps for the $T^2$ components and components $\Sigma\subset M$
without sutures,
correspond to sutures.
\ee

\section{Main Theorem}

\begin{thm}\label{maintheorem} Let $(M,\gamma)$ be an irreducible
sutured manifold
with  annular sutures, and let
$(M,\Gamma)$ be the associated convex structure. The following are equivalent.
\begin{enumerate}
\item[\rm(1)] $(M,\gamma)$ is taut.
\item[\rm(2)] $(M,\gamma)$ carries a taut foliation.
\item[\rm(3)] $(M,\Gamma)$ carries a universally tight contact structure.
\item[\rm(4)] $(M,\Gamma)$ carries a tight contact structure.
\end{enumerate}
\end{thm}

\proof
Without loss of generality we assume $M$ is connected.

(1)$\Rightarrow$(2) is Gabai's theorem \cite{Ga}.  Gabai's theorem does
not require the
assumption that $(M, \gamma)$ have annular sutures.

(2)$\Rightarrow$(1)
by Thurston
\cite{Th} does not require this assumption either.

(2)$\Rightarrow$(3) is due to Eliashberg and Thurston \cite{ET} in the closed case.
That their work can be applied in this context is the content of Theorem \ref{C}.

(3)$\Rightarrow$(4) is immediate.

(4)$\Rightarrow$(1) follows from Theorem \ref{B}.  The assumption that
$\bdry M \ne \emptyset$ is crucial here.  For by Bennequin \cite{Be83} $S^3$ has a
tight contact structure, but by Novikov \cite{No} it has no taut foliation.   Also the
irreducibility of $M$ is necessary, since connect summing preserve tightness (\cite{ML},\cite{Co97}), whereas
the universal cover of a taut foliation is $\R^3$. \qed

\subsection{Confoliations}

In this section we will prove the following theorem:

\begin{thm} \label{C}   Let $\xi$ be a (finite depth) taut foliation which is carried by a sutured
manifold
$(M,\gamma)$ with annular sutures.  Then there exists a modification of $\xi$ into a positive tight
contact structure $\xi_+$ such that ${\bdry M}$ is convex and $\Gamma_{\bdry M}=s(\gamma)$.
\end{thm}

Before we begin the proof, we recall several notions from the theory of confoliations \cite{ET}.
A {\it positive confoliation} $\xi$  is an oriented 2--plane field distribution on  $M$ given by a 1--form
$\alpha$ which satisfies $\alpha\wedge d\alpha\geq 0$.  The {\it contact part} of $\xi$ is
$H(\xi)=\{x\in M| \alpha\wedge d\alpha>0\}$.  For a subset $A\subset M$, the {\it saturation} $\hat A$ of
$A$ is the subset of $M$ which consists of points which can be
connected to a point in $A$ via a path which is everywhere tangent to $\xi$.
$\xi$ is said to be {\it transitive} if $\widehat{H(\xi)}=M$.

\proof The proof is almost identical to the perturbation result for closed manifolds due to
Eliashberg and Thurston \cite{ET}.  The difference is that we need to modify the boundary
carefully, and the modification $\xi_+$ is usually {\it not a perturbation} of $\xi$.
Since $\xi$ is carried by $(M,\gamma)$, $\bdry M$ is best thought of as a manifold with
corners, where $R_\pm=R_\pm(\gamma)$ are leaves of $\xi$ and the leaves of $\xi$ (and hence $R_\pm(\gamma)$)
are transverse to
$\gamma$.    In order to use symplectic filling techniques, we need to exercise a little care, and
extend $M$ and $\xi$ to an open manifold with {\it finite geometry at infinity}.

\vskip.15in
\noindent
{\bf Step 1}\qua  We first extend $\xi$ in two ways
to $M_1=M\cup (R_+\times [0,\infty))\cup
(R_-\times [0,\infty))$, where $R_+\times\{0\}=R_+$, $R_-\times\{0\}=R_-$,
$\bdry M_1=\gamma'$, and $\gamma'=\gamma\cup
(\bdry R_+\times [0,\infty))\cup(\bdry R_-\times [0,\infty))$ is smooth.
The first extension is to a foliation (still called $\xi$) and the second is to a positive
confoliation $\xi'$ which is contact on $R_\pm\times (0,\infty)$.
Let $t$ be the coordinate in the $[0,\infty)$--direction
for $R_+\times [0,\infty)$.
The extension to a foliation $\xi$ on $M_1$ is easy  -- on $R_\pm\times (0,\infty)$, simply
take $\ker dt$. We now construct $\xi'$.

\begin{lemma}
If $R_+$ has nonempty boundary, then there exists a 1--form $\beta$ on $R_+$ with $d\beta>0$, whose
singular foliation given by $\ker \beta$ has isolated singularities and no closed orbits,
and whose flow is transverse to $\bdry R_+$.
\end{lemma}

\proof  Start with a singular foliation $\cal{F}$ on $R_+$ which satisfies the
following:
\be
\item $\mathcal{F}$ is Morse--Smale and has no closed orbits,
\item The singular set consist of elliptic points (sources) and hyperbolic points.
\item $\mathcal{F}$ is oriented, and for one choice of orientation the flow is transverse to and exits from
$\bdry R_+$.
\ee
For example, a gradient-like vector field would do.
Next, modify $\mathcal{F}$ near each of the singular points so that $\mathcal{F}$ is given by
$\beta_0=ydx-xdy$ near an elliptic point and $\beta_0=ydx+2xdy$ near a hyperbolic point.  Therefore,
we have $\mathcal{F}$ given by $\beta_0$ which satisfies $d\beta_0>0$ near the singular points.
Now, let $\beta=f\beta_0$, where $f$ is a positive function with $df(X)>>0$, and $X$ is an oriented
vector field for $\mathcal{F}$ (nonzero away from the singular points).
Since $d\beta=df\wedge\beta_0+fd\beta_0$, $df(X)>>0$ guarantees that $d\beta>0$. \endproof

Choose a 1--form $\beta$ on $R_+$ as in the lemma.  Consider the 1--form $\alpha'=dt+f(t)\beta$
on $R_+\times[0,\infty)$, where $f(0)=0$, $f(t)=1$ for $t\geq 1$, and $f(t)>0$ for $t>0$.
$\alpha'\wedge d\alpha'=f(t) dt\wedge d\beta>0$
on $R_+\times(0,\infty)$, since $d\beta>0$.
Therefore, $\alpha'$ gives rise to an extension
of $\xi'$ to a positive confoliation on $M_1$.     The construction is similar on $R_-\times [0,\infty)$.
$\xi'$ is foliated on $M$ and contact on $M_1\backslash M$.

\vskip.15in
\noindent
{\bf Step 2}\qua Next extend $\xi$ to a foliation and $\xi'$ to a positive confoliation
on $M_2=M_1\cup (\gamma'\times [0,\infty))$.
Denote $\gamma'\times\{0\}=\gamma'$ and assign coordinates $(\theta,y,z)$ to
$\gamma'\times[0,\infty)=S^1\times \R\times[0,\infty)$ by setting
$y=\pm (t+1)$ on $\R_\pm\times[0,\infty)$ and $\gamma=S^1\times [-1,1]\subset \gamma'$.
Since $\xi|_\gamma=\xi'|_\gamma$ has no Reeb components, we may assume that
${\bdry \over \bdry y}\pitchfork \xi|_\gamma$.  This means that, on $\gamma'\times\{0\}$,
we can take the characteristic foliation for $\xi$ to be given by a 1--form
$\alpha=  dy-g(\theta,y,0)d\theta$, where $g=0$ if $y\geq 1$ or $y\leq -1$.   We extend $\alpha$ to a
foliated 1--form on $\gamma'\times I$ by taking $\alpha=dy-g(\theta,y,0)d\theta$.
Next, on $\gamma'\times \{0\}$, the characteristic foliation of $\xi'$
is given by the 1--form $\alpha'=dy-h(\theta,y,0)d\theta$, where $h<0$ for $y>1$ or $y<-1$,
and $h$ is independent of $y$ for large positive or large negative $y$.
Extend $\alpha'$ to a
positive confoliated 1--form  on $\gamma'\times [0,\infty)$ by taking $h$ with ${\bdry h\over \bdry z}<0$
and $\lim_{z\rightarrow \infty}h(\theta,y,z)=C$, where $C$ is a fixed large negative number.
Therefore, we have a confoliation $\xi'$ on $M_2$ whose contact part is $M_2\backslash M$.

Notice that if we took  $M\cup (R_\pm\times[0,1])\cup (\gamma''\times [0,n])$, $n$ large, where
$\gamma''=(\bdry R_+\times [0,1])\cup (\bdry R_-\times[0,1])\cup \gamma$, then we can round
the corners to obtain a manifold with boundary $M_3$ (isotopic to $M$ if we ignore corners).  The
characteristic foliation on $\bdry M_3$ is Morse--Smale, and $\Gamma_{\bdry M_3}$ is isotopic to
$s(\gamma)$.

\vskip.15in
\noindent
{\bf Step 3}\qua In this step we modify $\xi'$ on $M_2$
(fixing $\xi'$ on $M_2\backslash N(M)$, where $N(M)$ is a small neighborhood of $M$) to
obtain $\xi_+$ which is contact on all of $M_2$.
This step follows directly from Eliashberg and Thurston's argument \cite{ET}.
We list the relevant results:

\begin{prop}
Any $C^2$--confoliation can be $C^0$--approximated by a $C^1$--smooth transitive confoliation.
\end{prop}

\begin{prop}
Any $C^k$--smooth transitive positive confoliation, $k\geq 1$, admits a $C^k$--close approximation
by a positive contact structure.
\end{prop}

It is easy to see that the propositions hold while fixing $\xi'$ on $M_2\backslash N(M)$.  Therefore, we obtain
$\xi_+$ which is a positive contact structure and agrees with $\xi'$ `at infinity'.

\vskip.15in
\noindent
{\bf Step 4}\qua We prove that $(M_2,\xi_+)$ is symplectically semi-fillable.    We will construct a dominating
2--form $\omega$ for $\xi_+$ (ie, a closed 2--form for which $\omega|_{\xi_+}>0$ everywhere).

First recall the construction of a dominating 2--form $\omega$ on $M$ for the foliation
$\xi$.  Since
the foliation $\xi$ is taut, through each point there exists a closed
transversal or a transversal arc with endpoints on $R_+$ and $R_-$.  Let
$\delta_p$ be a transversal through the point $p$ and $N_p$ be a tubular
neighborhood of $\delta_p$.  Then $N_p$ is foliated by an interval's worth or
$S^1$'s worth of disks, and we have a projection $\pi_p\co  N_p\rightarrow D_p$,
where $D_p$ is a disk.  Let $\omega_p$ be the closed 2--form $\pi_p^*(f_pA_p)$,
where $A_p$ is an area form on $D_p$ and $f_p$ is a nonnegative function on $D_p$
with support inside $D_p$ and such that $f_p(\pi_p(p))>0$.  We may cover $M$ by
$N_p$ so that  $\bigcup_p supp (\omega_p)=M$, and take a finite subcover.
We would then take the dominating 2--form to be $\omega=\sum\omega_p$ (finite sum).
Note that these $\omega_p$ are additive.

For our purposes, we need to control this construction more carefully.
Let $M'=M\cup(R_+\times[0,\varepsilon])\cup (R_-\times[0,\varepsilon])$.   Extend the
transversal arcs $\delta_p$ ending at $R_\pm$ on $M$   so that on
$R_\pm\times[0,\varepsilon]$ they restrict to $\{pt\}\times [0,\varepsilon]$,
and choose $N_p$ so that $N_p\cap (R_+\times[0,\varepsilon])=
D_p\times [0,\varepsilon]$ (same for $R_-$).
Therefore, on $R_\pm\times[0,\varepsilon]$ we would have $\omega_p=\pi^*(g_pB_p)$,    where
$\pi\co R_\pm\times[0,\varepsilon]\rightarrow R_\pm$,
$B_p$ is some area form on $R_\pm$ and $g_p$ is a nonnegative function.
$\omega$ would then  have the property that
$\omega=\pi^*(A)$, where $A$ is some area form for $R_\pm$.  Therefore we can
extend $\omega$ to $M_1$ so that
$\omega=\pi^*(A)$, where $\pi\co R_\pm\times [0,\infty)\rightarrow R_\pm$ is the first projection and
$A$ is an area form for $R_\pm$.
We can further extend it to $M_2$ so that
$\omega=dzd\theta$ on $\gamma'\times [\varepsilon,\infty)$.
Extending in the $M_2$--direction is easy
if we took care to choose (1) $\delta_p$ to be arcs with $\theta=const.$ and $z=0$, if $p\subset \gamma'$,
and (2) $N_p\subset M'\cup (\gamma'\times[0,\varepsilon])$.  This means we can simply
add the form $f(z)dzd\theta$, where $f(z)=1$ for $z\geq \varepsilon$, $f(0)=0$, and $f(z)> 0$
for $z> 0$.
By our construction of
$\xi_+$, the closed 2--form $\omega$ satisfies
$\omega|_{\xi_+}>0$ as well as $\omega|_{\xi}>0$.

Define a closed 2--form
$\widetilde \omega= \omega+d(s\alpha)$  on $M_2\times [-\varepsilon,\varepsilon]$,
where $s$ is the variable for $[-\varepsilon,\varepsilon]$, $\alpha$ is a nowhere zero
1--form whose kernel is $\xi$, and $\varepsilon>0$ is small enough.   Since we can obtain $\xi_+$ positive and
$\xi_-$ negative (similarly), $(M_2,\xi_+)$ is symplectically semi-fillable and dominated by $\widetilde\omega$.
We have the following symplectic semi-filling result:

\begin{thm} [Gromov--Eliashberg]   \label{filling}
Let $(X,\widetilde\omega)$ be a (not necessarily compact) symplectic 4--manifold with
contact boundary $(M,\xi)$ which satisfies $\widetilde\omega|_\xi>0$.  Assume there exists a calibrated
almost complex structure $J$ on $M$ which preserves $\xi$, and a corresponding
Riemannian metric $g$ which has {\em finite geometry at infinity}, ie,
\be
\item $g$ is complete,
\item the sectional curvature of $g$ is bounded above, and
\item the injectivity radius of $g$ is bounded below by some $\varepsilon>0$.
\ee
Then $(M,\xi)$ is a tight contact manifold.
\end{thm}

By our construction, $M_2\times[-\varepsilon,\varepsilon]$ has {\it finite geometry at infinity}.
Now pass to the universal cover of $M_2\times [-\varepsilon,\varepsilon]$, which also has
finite geometry at infinity. Theorem \ref{filling} implies that $\xi_+$ is universally tight.
Hence so is $\xi_+$ restricted to $M_3$.\endproof

\noindent
{\bf Remark}\qua It is possible to prove that if $(M,\gamma)$ is taut, then $(M,\Gamma)$ carries a universally
tight contact structure without resorting to symplectic filling.
Instead we may use a convex decomposition which matches Gabai's sutured manifold decomposition,
and prove a gluing theorem for tight contact structures.  This will be carried out in \cite{HKM}, using ideas
in \cite{H2}.

\subsection{Proof of (4)$\Rightarrow$(1)}

\begin{thm}\label{B} If $(M,\Gamma)$ carries a tight contact structure
then $(M,\Gamma)$ is taut.
\end{thm}

\proof  Let us assume instead that there exists a surface $T\subset
M$ such that
\be
\item $[T]=[R_+(\Gamma)]=[R_-(\Gamma)]\in H_2(M,\Gamma)$.
\item $x(T)< x(R_+(\Gamma))$.
\ee
The proof will follow from a sequence of lemmas and a calculation in the end.

\begin{lemma}\label{kubota} It is possible to modify $T$ so that $T$ 
satisfies {\rm(1)},
{\rm(2)} as well as
\be
\item[\rm(3)] $\bdry T=\Gamma$.
\ee
\end{lemma}

\proof  Let $\Gamma_0$ be a connected component of $\Gamma$, and consider all
the `sheets' $T_1,\cdots,$ $T_m$ of $T\cap N(\Gamma_0)$, where
$N(\Gamma_0)$ is a small neighborhood of
$\Gamma_0$. Since $[R_+(\Gamma)]\mapsto [\Gamma]$ under the boundary map
$H_2(M,\Gamma)\stackrel{\bdry}{\rightarrow}H_1(\Gamma)$, if $m>1$, then
there must exist two consecutive
sheets $T_i$ and $T_{i+1}$ which are oppositely oriented.  In this
case, we may surger $T$
by gluing $T_i$ and $T_{i+1}$ along $\Gamma_0$, rounding,
and pushing the two sheets off of $\Gamma_0$.
In this fashion we may reduce $m$ until it eventually becomes $1$. \qed

\begin{lemma}\label{x=chi} In addition, we may take $T$ to satisfy
\be
\item[\rm(4)] $x(T)=-\chi(T)$.
\ee
\end{lemma}

\proof This is asking that $T$ have no disk or sphere components, which are
the ones that contribute positively to the Euler characteristic but do not
contribute to the Thurston norm. The irreducibility of $M$ assures us that
every $S^2$ bounds a 3--ball, and
can  be removed from $T$ without affecting homology.
We claim that there can be no disks $D$ with $\delta=\bdry D$ which is a component of
$\Gamma$, unless $(M,\Gamma)=(B^3,S^1)$.
If there is such a disk $D$, then take a curve $\delta'\subset \bdry M$ parallel to
$\delta$ which has no intersections with $\Gamma$.  Use the Legendrian Realization
Principle to realize $\delta'$ as a Legendrian curve with $t(\delta', \bdry M)=0$.
$\delta'$ will then bound a disk $D'$ with $t(\delta',D')=0$.  This is an equivalent
definition of the existence of an overtwisted disk.
If $(M,\Gamma)=(B^3,S^1)$, Theorem \ref{B} is immediate.\qed

\begin{lemma}\label{defineW}  In addition, $T$  may be modified so 
that
\be
\item[\rm(5)] $W$, the
union of components $M\backslash T$ which intersect $R_+(\Gamma)$, satisfies
$\partial W  = R_+(\Gamma) \cup T_-$.
\ee
\end{lemma}

Here, if $M'=M\backslash T$, then we define $T_+$, $T_-$ to be copies of $T$ contained
in $M'$, where the orientation induced by $T$ points out of and into $M'$ (respectively).

\proof Define the function $\phi\co  M\backslash T\rightarrow \Z$ as
follows. Assign $\phi(M_0)=0$, where $M_0$ is some connected component of
$M\backslash T$ which borders $R_-(\Gamma)$.  For another connected
component $M_i$, take an arc $\alpha$ which starts in $M_0$ and ends in
$M_i$, and define $\phi(M_i)=[\alpha\cap T]$. This number is independent of
$\alpha$, for if $\alpha'$ is another  curve with the same endpoints, then
$[(\alpha-\alpha')\cap T]=[(\alpha-\alpha')\cap  R_+(\Gamma)]=0$. Note that
all the components $M_i$ which border $\partial M$ have $\phi(M_i)$ equal to
either $0$ or $1$. If $\phi(M\backslash T)\not = \{0,1\}$, then choose $M_i$
with $\phi(M_i)$ extremal. $\partial M_i$ will not intersect $\partial M$ and
will consist of  components of $T$.  Since these components bound $M_i$, we
may throw them away without increasing $x(T)$.  Thus we may assume
$\phi(M\backslash T)=\{0,1\}$.

Let $M_1$ be a component of $M\backslash T$ which intersects 
$R_+(\Gamma)$ and let
$\alpha$ be an arc which starts in $M_0$ and ends in $M_1 \cap 
R_+(\Gamma)$.  Then
$\phi(M_1) = [\alpha \cap T] = [\alpha \cap R_+(\Gamma)] = 1$.  $M_1 \cap
R_-(\Gamma) = \emptyset$, for otherwise there exists an arc $\beta$ connecting
points of $R_+(\Gamma)$ and $R_-(\Gamma)$ which doesn't intersect 
$T$.  Also $M_1
\cap T_+ = \emptyset$, for crossing $T_+$ increases $\phi$, and $\phi$
already takes
its maximum value on $M_1$.  It follows that $\partial W_1 \subset 
R_+(\Gamma) \cup T_-$.

Conversely, suppose that $M_2$ is a component of $M\backslash T$ 
which intersects
$T_-$.  Since crossing $T_-$ decreases $\phi$, it follows that 
$\phi(M_2)= 1$.  $M_2
\cap R_-(\Gamma) = \emptyset$; otherwise following an arc from $R_-(\Gamma)$ to
$R_+(\Gamma)$ would increase the value of $\phi$ by 1.  Also $M_2 \cap T_+ =
\emptyset$ since crossing $T_+$ increases $\phi$.  For $M_2$ to be included in
$W$, we require that $M_2$ intersect $R_+(\Gamma)$.  If this is not 
the case, then
$\partial M_2 \subset T_-$, and this component of $T$ can be 
eliminated from $T$.
\qed

\begin{lemma} \label{goku}
There exists an isotopy $\phi_t\co T\rightarrow M$, $t\in[0,1]$, such that
$\phi_0(T)=T$,  $S\stackrel{\rm def}{=}\phi_1(T)$ is a convex surface, and
$\phi_t(\partial T)$, $t\in[0,1]$, is contained in an annulus $N(\Gamma)\subset \bdry M$
which contains $\Gamma$.
\end{lemma}

\proof By Lemma~\ref{kubota}, $\partial T = \Gamma$.  Perturb $T$ so that each
component of $\partial T$ is transverse to and non-trivially intersects $\Gamma$.  By
the Legendrian Realization Principle we may assume $\partial T$ is a union of
Legendrian curves.  By Theorem~\ref{existence}, $T$ may be
isotoped to a convex surface.
\endproof

\noindent \textbf{Completion of the proof of Theorem~\ref{B}}\qua Let
$(M',\Gamma')$ denote $(M, \Gamma)$ split along $(S, \sigma)$, where
$S$ is as in Lemma \ref{goku} and $\sigma$ is its dividing set.  Also
let $W$ be as in Lemma~\ref{defineW}. Recall $\partial W = R_+(\Gamma) \cup T_-$.
By our choice of $S$, $M\backslash S \cong M\backslash T$, and
we denote the components of  $M\backslash S$ which correspond to $W$ by
$\bar W$. The convex structure on $\bar W$ is denoted $(\bar W, \bar
\Gamma, R_-(\bar \Gamma), R_+(\bar \Gamma))$.

We must show how $R_+(\bar \Gamma)$ is related to $R_+(\Gamma)$.  Let
$N(\Gamma)$ be a
regular neighborhood of $\Gamma$ in $\partial M$ which contains the isotopy of
$\partial T$ to $\partial S$.  Let $R_+$ be the closure of $R_+(\Gamma)\backslash
N(\Gamma)$.  It follows that $R_+$ is contained in the interior of $R_+(\bar
\Gamma)$.  It follows that there exist subsurfaces $A$ and
$B$ of $\partial \bar W$ which intersect along circles such that
$$R_+(\bar\Gamma) = R_+ \cup A$$
$$R_-(\bar\Gamma) = B$$
$$A \cup B \cong S_- \cong T_-.$$
By Corollary~\ref{A1}, $\chi(R_+(\bar\Gamma)) = \chi(R_-(\bar \Gamma))$.  An
argument similar to that of the proof of
Lemma \ref{x=chi} gives $\chi(A) \le 0$; thus it follows that
$$\chi(R_+) =  \chi(B) - \chi(A) \ge \chi(B) + \chi(A) = \chi(T_-).$$ 
Since $R_+
\cong R_+(\Gamma)$ and $T_- \cong T$, it follows that $x(R_+(\Gamma)) 
\le x(T)$. \endproof

{\bf Acknowledgements}\qua  We thank the referee for helpful comments on improving the
exposition.


\end{document}